\documentclass[12pt]{amsart}
\usepackage[english]{babel}

\usepackage{amssymb,amsmath,amscd,amsthm}

\def\{{\protect\lbrace}
\def\}{\protect\rbrace}

\newcommand{\Aut}{\operatorname{Aut}}

\begin{document}
\begin{center} 
\textbf{Rings with Polynomial\\ Identity and Centrally Essential Rings}
\end{center} 
\hfill {\sf V.T. Markov}

\hfill Lomonosov Moscow State University

\hfill e-mail: vtmarkov@yandex.ru

\hfill {\sf A.A. Tuganbaev}

\hfill National Research University "MPEI"

\hfill Lomonosov Moscow State University

\hfill e-mail: tuganbaev@gmail.com

{\bf Abstract.} It is proved that for any prime integer $p$ and each field $F$ of characteristic $p$, there exists a centrally essential $F$-algebra which is not a PI-ring and is not algebraic over its center.

V.T.Markov is supported by the Russian Foundation for Basic Research, project 17-01-00895-A. A.A.Tuganbaev is supported by Russian Scientific Foundation, project 16-11-10013.


\textbf{MSC2010 datebase 16R99; 16D10}

\section{\bf Introduction}\label{section1}
All considered rings are associative and contain the non-zero identity element. 

\textbf{1.1. Centrally essential rings.} A ring $R$ with center $C=C(R)$ is said to be \textsf{centrally essential}\footnote{Centrally essential rings are studied, for example, in \cite{MT18}.} if the module $R_C$ is an essential extension of the module $C_C$, i.e., for any non-zero element $r\in R$, there exist non-zero central elements $c,d\in C$ with $rc=d$. 

It is clear that all commutative rings are centrally essential. If $Z_2$ is the field of order 2 and $Q_8$ is the quaternion group of order 8, then the group ring $Z_2[Q_8]$ is an example of a non-commutative, centrally essential, finite ring \cite{MT18}.

\textbf{1.2. Rings with polynomial identity.} Let $X$ be a countable set and $F=\mathbb{Z}<X>$ a free ring with the set of free generators $X$. A \textsf{classical identity} in the sense of Rowen is an identity with integral coefficients, i.e., the element of the free ring $F$ which is contained in the kernel of any homomorphism of the ring $F$ in ring $R$. 
A classical identity is called a \textsf{polynomial identity} if it is multilinear and has $1$ as one of its coefficients; a ring with polynomial identity is called a 
\textsf{PI ring}\footnote{See \cite[Definitions 1.1.12, 1.1.17]{Row80}}. 

The main result of this paper is Theorem 1.3.
 
\textbf{1.3. Theorem.} For any prime integer $p$ and each field $F$ of characteristic $p$, there exists a centrally essential $F$-algebra which is not a PI ring and is not algebraic over its center.

\textbf{1.4. Rings which are algebraic or integral over their centers.}\\
Let $R$ be an arbitrary ring with center $C$. An element $r\in R$ is said to be \textsf{algebraic} (resp., \textsf{integral}) over the center if, for some $n\in \mathbb{N}$, there exist elements $c_0,\ldots c_n\in C$ such that $c_n$ is a non-zero-divisor in $R$ (resp., an invertible  element in $R$) and
\begin{equation}\label{algeq}
c_nr^n+c_{n-1}r^{n-1}+\ldots+c_{1}r+c_0=0.
\end{equation} 
We denote by $n_{1}(r)$ (resp., $n_{2}(r)$) the least integer $n$ which satisfies this condition. A ring $R$ is said to be \textsf{algebraic} (resp., \textsf{integral}) over its center if any element $r\in R$ is algebraic (resp., \textsf{integral}) over its center. We set $m_{1}(R)=\max\{n_{1}(r)\,|\,r\in R\}$ and $m_{2}(R)=\max\{n_{2}(r)\,|\,r\in R\}$; it is possible that $m_{1}(R)=\infty$, $m_{2}(R)=\infty$. 

Finite rings and finite-dimensional algebras over fields are examples of rings $R$ such that $m_1(R)=m_2(R)<\infty$. 

\textsf{We give an example of the ring which is algebraic and is not integral over its center.} Let 
$$R=\left\lbrace \begin{pmatrix}
a&b\\0&z
\end{pmatrix}\colon a,b\in \mathbb{Q},\, z\in\mathbb{Z}
\right\rbrace,
$$
It is clear that the center of the ring $R$ is of the form $\mathbb{Z} E$, where $E$ is the identity matrix.

We note that the ring $R$ is not integral over its center. Indeed, since $\mathbb{Q}$ is a homomorphic image of the ring $R$, the integrity of $R$ over its center would imply that $\mathbb{Q}$ is integral over $\mathbb{Z}$ which is false.

On the other hand, if $r\in R$, then $nr\in M_2(\mathbb{Z})$ for some $n\in\mathbb{N}$, and the matrix ring $M_2(\mathbb{Z})$ is a finitely generated $\mathbb{Z}$-module and hence it is integral over $\mathbb{Z}$.

We note that the classes of centrally essential rings, PI rings and rings which are algebraic (or integral) over its center, properly contain the class of all commutative rings.

\section{\bf The Proof of Theorem 1.3}\label{section2}

The proof of Theorem 1.3 uses the following two familiar results.

\textbf{2.1. Theorem.} Passman, 1977 \cite[Theorem 5.3.9(ii)]{Passman}.
Let $F$ be a field of characteristic $p>0$. If the group algebra $FG$ satisfies a polynomial identity of degree $d$, then there exists a subgroup $H$ of the group $G$ such that $[G:H]\cdot|H'|<g(d)$, where $g(d)$ is some fixed function of the integer $d$.

\textbf{2.2. Theorem.} Markov--Tuganbaev, 2018 \cite[Theorem 1.1.2]{MT18}.
Let $F$ be a field of characteristic $ p > 0$. If $G$ is a finite \hbox{$p$-group} of nilpotence class 2, then the group algebra
$FG$ is a centrally essential ring.

We fix a prime integer $p$ and a field $F$ of characteristic $p$. 
We denote by $Z(G)$ the center of the group~$G$.

\textbf{2.3. Lemma.} There exists a series of finite $p$-groups $G(n)$, $n\in \mathbb{N}$, such that\\
1) for any $n\in \mathbb{N}$, the group algebra $FG(n)$ is a centrally essential ring;\\
2) for any $d\in \mathbb{N}$ there exists $n=n(d)$ such that
the ring $FG(n)$ does not satisfy a polynomial identity of degree $d$.

\begin{proof}
For any positive integer $n$, we construct a group $G=G(n)$ as follows. 
Let $A=\langle a\rangle$, $B=\langle b\rangle$, $C=\langle c\rangle$ be three cyclic groups such that 
$|A|=|B|=|C|=p^n$. We consider the automorphism $\alpha\in \Aut(B\times C)$ defined on the generators by the relations $\alpha(b)=bc$ and $\alpha(c)=c$. It is clear that $\alpha^n$ is the identity automorphism; 
therefore, we have a homomorphism $\varphi\colon A\rightarrow\Aut(B\times C)$ with $\varphi(a)=\alpha$. This 
homomorphism corresponds to the semidirect product $G=(B\times C)\ltimes A$ which can be 
considered as the group generated by the elements $a,b,c$ which satisfy the relations
$a^{p^n}=a^{p^n}=c^{p^n}=1$, $bc=cb$, $ac=ca$ and $aba^{-1}=bc$. 
It follows from these relations that $c\in Z(G)$. 
It is directly verified that for any integers $x,y,z, x',y',z'$, we have
\begin{equation}
\label{comm}
\begin{array}{rcl}[b^yc^za^x,b^{y'}c^{z'}a^{x'}]&=
&b^ya^xb^{y'}a^{x'}a^{-x}b^{-y}a^{-x'}b^{-y'}\\
&=&b^y(a^xb^{y'}a^{-x})(a^{x'}b^{-y}a^{-x'})b^{-y'}\\
&=&b^y(b^{y'}c^{xy'})(b^{-y}c^{-yx'})b^{-y'}=c^{xy'-yx'}.\end{array}
\end{equation}
Thus, $Z(G)=G'=\langle c\rangle$ and $G$ is a group of nilpotence class 2, so the first assertion follows from Theorem 2.1.

Let $H$ be an arbitrary subgroup of the group $G$. Then

\begin{equation}
\label{ineq}
[G:H]\cdot|H'|\geq p^n.
\end{equation}

We note that $[G:HZ(G)]\leq [G:H]$ and $(HZ(G))'=H'$; consequently, it is sufficient to prove the inequality \eqref{ineq} in the case, where $H\supseteq Z(G)$. We set $\bar G=G/Z(G)$ and denote by $\bar a, \bar b, \bar H$ the images of $a,b,H$ under the canonical homomorphism $G$ onto the group $\bar G$. We also set $\bar B=\langle \bar b\rangle$. We have $[G:H]=[\bar G:\bar H]$. It follows from the standard isomorphism 
$(\bar H\bar B)/\bar B\cong \bar H/(\bar H\cap \bar B)$ that $\bar H/(\bar H\cap\bar B)$ is a cyclic group which is to some subgroup of the group $\langle \bar a\rangle$. The group $\bar H\cap\bar B$ is also cyclic; consequently, the group $\bar H$ is generated by two elements
of the form $\bar b^{p^m}$ and $\bar a^{p^k}\bar b^l$ for some non-negative integers $k,l,m$. Thus,
$$[\bar G:\bar H]=[\bar G:\bar H\bar B][\bar H\bar B:\bar H]=[\langle \bar a\rangle:
\langle \bar a^{p^k}\rangle] [\langle \bar b\rangle: \langle \bar b^{p^m} \rangle] =p^kp^m=p^{m+k}.$$

If $m+k\geq n$, then the inequality \eqref{ineq} holds. If $m+k<n$, then, by \eqref{comm} and the property that the elements $a^{p^k}b^l$ and $b^{p^m}$ are contained in the subgroup $H$, we have that $[a^{p^k}b^l,b^{p^m}]=c^{p^{m+k}}\in H'$; therefore, $|H'|\geq |\langle c^{p^{m+k}}\rangle|=p^{n-m-k}$ and we have
$$
[G:H]\cdot|H'|\geq p^{m+k}\cdot p^{n-m-k}=p^n,
$$
i.e., \eqref{ineq} also holds in this case.

Hence  the second assertion follows from Theorem 2.2.
\end{proof}

Now we finish the proof of Theorem 1.3. It is sufficient to take the direct product of the group algebras $FG(n)$, $n\in \mathbb{N}$, satisfying the conditions of Lemma 2.3,  as the ring $R$. We note that the direct product of any set of rings is centrally essential if and only if every factor is a centrally essential ring. Therefore, the ring $R$ is centrally essential. However, if the algebra $R$ satisfies some polynomial identity of degree $d$, then every ring $FG(n)$ satisfies this polynomial identity, contrary to the second assertion of Lemma 2.3.

It remains to prove that the constructed ring is not algebraic over its center.

It is well known (e.g., see \cite[Proposition 1.1.37]{Row80} or \cite[Lemma 5.2.6]{ZhSSS}) that if $m_1(R)=m<\infty$, then $R$ satisfies a polynomial identity of degree $\displaystyle{d(m)=\frac{m(m+1)}{2}+m}$.

We note that for any $m\in \mathbb{N}$, there exists an integer $n_m$ such that $R(n_m)$ does not satisfy any polynomial identity of degree $d(m)$; moreover, we can choose integers $n_1,n_2,\ldots$ such that these integers form an ascending sequence. By the definition of $d(m)$, there exists an element $r'_m\in FG(n_m)$ which does not satisfy any relation of the form \eqref{algeq} of degree $m$. Now we consider the element
$r=\prod_{n=1}^\infty r_n\in \prod_{n=1}^\infty FG(n)$,
where $r_n\in FG(n)$, $r_n=r'_m$ if $n=n_m$ for some $m\in\mathbb{N}$ and $r_n=0$, otherwise. It is clear that if $r$ satisfies some relation
of the form \eqref{algeq} of degree $m$, then every element $r_n$ satisfies a relation of the same degree; this is  impossible by the choice of the element 
$r'_m$.~\hfill$\square$

\textbf{2.4. Remark.}
The reviewer suggested another example of a centrally essential ring which is not a PI ring, namely the group 
algebra $FG$ where $G$ is the direct sum\footnote{The direct sum of some 
groups $G(n)$ is the subgroup of 
$\prod_{n\in\mathbb{N}}G(n)$ consisting of such elements $\prod_{n\in\mathbb{N}}g_n$ that only finite number of elements $g_n$ are non-identical.} of all the groups $G(n)$ 
constructed in the proof of Lemma 2.3. It is easy to check 
that this algebra is a centrally essential ring. It is also clear that it is not a PI ring since every ring $FG(n)$ is its subring (as well as its homomorphic image). But this ring is evidently  integral over its center since is locally finite over $F$.

\textbf{2.5. Acknowledgment.} The authors are sincerely grateful to the reviewer for valuable comments and suggestions.

\end{document}